\documentclass[11pt,twoside]{amsart}

\usepackage{amsmath, amsthm, amssymb, mathrsfs}
\usepackage{array}
\usepackage{hyperref}
\usepackage[normalem]{ulem}
\usepackage{color}
\usepackage{verbatim,tikz}

\input xy
\xyoption{all}

\numberwithin{equation}{section}
\newtheorem{theorem}{Theorem}[section]
\newtheorem{lemma}[theorem]{Lemma}
\newtheorem{proposition}[theorem]{Proposition}
\newtheorem{corollary}[theorem]{Corollary}

\theoremstyle{definition}
\newtheorem{definition}[theorem]{Definition}

\newtheorem{example}[theorem]{Example}

\DeclareMathOperator\cone{Cone}

\DeclareMathOperator\conv{conv}

\DeclareMathOperator\rank{rank}

\DeclareMathOperator\rr{\mathbb{R}}

\DeclareMathOperator\zz{\mathbb{Z}}

\DeclareMathOperator\fff{\mathbb{F}} 

\DeclareMathOperator\ppp{\mathcal{P}}

\DeclareMathOperator\lb{\langle}
\DeclareMathOperator\rb{\rangle}

\newcommand{\nr}{N_{\rr}}

\newcommand{\s}{\sigma}
\renewcommand{\ss}{\Sigma}

\newcommand{\Ga}{\Gamma}
\newcommand{\ib}{\overline{I}}
\newcommand{\llm}{\{1,\ldots,m\}}
\newcommand{\lln}{\{1,\ldots,n\}}

\newlength{\bracewidth}
\newcommand{\myunderbrace}[2]{\settowidth{\bracewidth}{$#1$}#1\hspace*{-1\bracewidth}\smash{\underbrace{\makebox{\phantom{$#1$}}}_{#2}}}

\makeatletter
\renewcommand*\env@matrix[1][*\c@MaxMatrixCols c]{%
  \hskip -\arraycolsep
  \let\@ifnextchar\new@ifnextchar
  \array{#1}}
\makeatother

\begin{document}

\title{Primitive Collections for Klyachko Varieties}

\author{Nicole Lemire, Jason Palombaro}

\begin{abstract}
Klyachko~\cite{Kly83} introduced a family of $S_m\times S_n$-invariant toric Gorenstein Fano varieties, $X_{mn}$, with connections to representation theory, algebraic groups and graph theory.  He focused on the case in which $m$ and $n$ are coprime, for which $X_{mn}$ are smooth, and hence Fano.
Klyachko and Voskresenskii~\cite{VK84} later proved that all forms of $X_{mn}$ are stably rational when $m,n$ are coprime, and that all forms of $X_{mn}$ are rational when $n$ is congruent to 1 modulo $m$. In the same paper, they classified centrally symmetric toric Fano varieties as products of projective lines and generalized del Pezzo varieties, $X_{2,2k+1}$.  For this reason, the special case of generalized del Pezzo varieties are frequently studied. For example, Ballard, Duncan, and McFaddin~\cite{BDM22} proved that the bounded derived category of  generalized del Pezzo varieties has an  full exceptional collection of line bundles.  An understanding of the maximal cones and primitive collections of the fan associated to $X_{2,2k+1}$ is  crucial in their construction.
In this paper, we focus on the Klyachko varieties $X_{mn}$ in the case in which $n$ is congruent to 1 modulo $m$.  We determine a description of the maximal cones and primitive collections for the associated fan.  We also determine a description of their cohomology ring, and some stability results on the $S_m\times S_n$-orbits of maximal cones and primitive collections for $X_{mn}$.
\end{abstract}

\maketitle

\section{Introduction}

Klyachko~\cite{Kly83} introduced a family of $S_m\times S_n$-invariant toric Gorenstein Fano varieties, $X_{mn}$, with connections to representation theory, algebraic groups and graph theory.  He focused on the case in which $m$ and $n$ were coprime, for which $X_{mn}$ were smooth, and hence Fano.
His motivation was to determine smooth projective models for a certain special class of non-split algebraic $k$-tori.

Soon afterwards, Klyachko and Voskresenskii~\cite{VK84} classified centrally symmetric toric Fano varieties, those whose anticanonical divisor is ample and whose corresponding anticanonical polytope is additionally symmetric.  
They showed that all centrally symmetric toric Fano varieties are products of projective lines and generalized del Pezzo varieties.
The generalized del Pezzo varieties correspond to the special case of Klyachko varieties with $V_{2k}=X_{2,2k+1}$.

Many properties of normal separated toric varieties can be determined by studying the underlying fan from which they are constructed. The primitive collections of the fan are of particular interest as they can be used to compute cohomology of the toric variety, 
construct the toric variety as a (almost) geometric quotient, 
and can play a key role in finding exceptional collections for the bounded derived category of the toric variety. They are closely tied to the maximal cones of the fan as they can be algorithmically determined if it is known which rays generate maximal cones: find the smallest possible sets that are not contained in any maximal cone. However, it is not always known which rays generate maximal cones, and when the fan is not fixed, such as when we are working with an infinite family of fans, we must resort to other measures for finding primitive collections.

The generalized del Pezzo varieties form an infinite family of varieties that are frequently studied because of their symmetry. The maximal cones of the fan associated to these varieties were found by Voskresenskii and Klyachko~\cite{VK84}[Theorem 5] and Ballard, Duncan, and McFaddin identified which subsets of rays form primitive collections from this result. These primitive collections were crucial in the proof of the main result of~\cite{BDM22}: a full exceptional collection of line bundles for the bounded derived category of coherent sheaves (\cite{BDM22}[Theorem 3.3]). Full exceptional collections give significant insight into the derived category, so it would be interesting to know whether this result can be generalized to a larger class of varieties.

A natural generalization of the generalized del Pezzo varieties is the class of Klyachko varieties $X_{mn}$.
For general $m,n$, $X_{mn}$ is a Gorenstein Fano toric variety with terminal singularities (\cite{CFST18}). 
These highly symmetric varieties were studied in~\cite{VK84} where the authors showed that $X_{mn}$ is smooth if and only if $m$ and $n$ are relatively prime.  If $m$ and $n$ are coprime, Voskresenskii and Klyachko proved that any  form of $X_{mn}$ is stably rational, and if $n=1 \mod m$ (or $m=1 \mod n$) then any form is rational. Since the case $n=1 \mod m$ has proven to be special and includes the generalized del Pezzo varieties, it makes for an ideal starting point for generalizing the results of~\cite{BDM22}.

The main result of this paper is a description of the maximal cones and primitive collections for the fans associated to the Klyachko varieties $X_{mn}$ for which $n=1 \mod m$. This is a significant first step toward finding exceptional collections of line bundles for these varieties and provides insight towards finding the maximal cones and primitive collections for the more general case in which $m$ and $n$ are only relatively prime.

Recall that root systems arise naturally in determining the maximal tori of algebraic groups.  In particular, the maximal torus of $\mathrm{PGL}_m$ has character lattice the root lattice of the root system $A_{m-1}$ and the maximal torus of $\mathrm{SL}_m$ has character lattice the weight lattice of the root system $A_{m-1}$.
The character lattice of $X_{mn}$ is the tensor product of the root lattice of the root system $A_{m-1}$ with the root lattice of the root system $A_{n-1}$.  Its cocharacter lattice $N$ can be identified with a quotient of $M_{mn}(\zz)$.  Let $e_{ij}\in N$ be the image of the natural basis element $E_{ij}$ of $M_{mn}(\zz)$. The Klyachko variety $X_{mn}$ is the projective toric variety corresponding to the  face fan of the reflexive polytope $P_{mn}$ determined as the convex hull of $\{e_{ij}:1\le i\le m, 1\le j\le n\}$.  The natural permutation action of $S_m\times S_n$ on $M_{mn}(\zz)$
induces an action on $P_{mn}$ and hence on $X_{mn}$.

The primitive collections for $X_{mn}$ in the case $n=km+1$ can be arranged into $m$ orbits  under the $S_m\times S_n$ action each with a standard `rectangular' representative for each `row' $1\le y\le m$:
$$\{e_{ij}: 1\le i\le y, 1\le j\le n-ky\}$$
A more precise statement is found in Proposition 3.8, along with its proof.

Proposition 3.8 also gives a graph theoretic description of the maximal cones in the fan of $X_{mn}.$
The vertices of the polytope $e_{ij}=\overline{e_i}\otimes \overline{f_j}$ determine edges in the complete bipartite graph $\Gamma_{mn}$ with vertex sets $\{\overline{e_i}:i=1,\dots,m\}$ and $\{\overline{f_j}:j=1,\dots,n\}$.
The maximal cones of $\Sigma_{mn}$ are determined by the spanning trees of $\Gamma_{mn}$ such that the degree of each vertex $\overline{e_i}$ is equal to $k+1=\frac{n+m-1}{m}$.
If $I$ is the edge set of such a spanning tree, then $\cone\{e_{ij}: (i,j)\not\in I\}$ is the corresponding maximal cone.

This paper begins with a detailed construction of the Klyachko varieties $X_{mn}$ and an explanation of the known geometric properties that they are known to satisfy from the literature~\cite{VK84,CFST18}.  We then introduce matrix representations which provides a way to write subsets of ray generators as a matrix which behaves nicely under the group action. 
Next we work towards the simultaneous characterization of the maximal cones and primitive collections for the fan $\ss_{mn}$, culminating in the result of Proposition 3.8.

After illustrating our results through the example corresponding to the smallest non-trivial case $X_{34}$, we prove some stability results on the orbits of maximal cones under the group action.
In particular, we prove that the orbits of primitive collections of $\Sigma_{m,km+1}$ under the action of $S_m\times S_{km+1}$ are in bijection with those of $\Sigma_{m,m+1}$ under the action of $S_m\times S_{m+1}$ for all $m\ge 2$, whereas the corresponding statement for maximal cones only holds for $m=2$ and $m=3$.  For $m\ge 3$, $k\ge m-2$, the orbits of maximal cones of $\Sigma_{m,km+1}$ under the action of $S_m\times S_{km+1}$ are in bijection with those of $\Sigma_{m,m(m-2)+1}$ under the action of $S_m\times S_{m(m-2)+1}$. For comparison, there is a unique orbit of maximal cones of $\Sigma_{2,2k+1}$ under the action of $S_2\times S_{2k+1}$.

As a last result, we determine a bijection of the orbits of maximal cones in the special case $\Sigma_{m,m+1}$ under the action of $S_m\times S_{m+1}$ with the set of unlabelled spanning trees on $m+1$ vertices.

The last section is an application of the main result: we use the found primitive collections to compute the cohomology ring for $X_{mn}$ in the case that $n=1 \mod m$.

In a subsequent paper, based partly on work in~\cite{Pal}, we apply these results to determine $S_m\times S_n$-invariant exceptional collections of line bundles for the Klyachko varieties, $X_{mn}$ with $n=1 \mod m$ following~\cite{BDM22}.

\section{The Klyachko Varieties}

\subsection{Construction}

For a finite group $G$, with subgroup $H$, let the transitive $G$-permutation lattice with stabilizer subgroup $H$ be $\zz[G/H]$.  We may identify the basis of $\zz[G/H]$ with the cosets $\{gH:g\in G/H\}$.
Let the augmentation map be $$\epsilon_{G/H}:\zz[G/H]\to \zz, gH\to 1.$$ 
Let $I_{G/H}=\ker(\epsilon_{G/H})$ be the kernel of the augmentation ideal and $J_{G/H}=I_{G/H}^*$ be its $\zz$-dual.
Then clearly:
$$0\to I_{G/H}\to \zz[G/H]\stackrel{\epsilon_{G/H}}{\to} \zz\to 0$$
$$0\to \zz\stackrel{\epsilon_{G/H}^*}\to \zz[G/H]\to J_{G/H}\to 0$$
are short exact sequences of $G$-lattices.
Note that $\epsilon_{G/H}^*$ sends 1 to the sum of the permutation basis in $\zz[G/H]$.

Recall that root systems arise naturally in determining the maximal tori of algebraic groups.  In particular, the maximal torus of $\mathrm{PGL}_m$ has character lattice the root lattice of the root system $A_{m-1}$ and the maximal torus of $\mathrm{SL}_m$ has character lattice the weight lattice of the root system $A_{m-1}$.
Note that the root lattice $\zz A_{m-1}$ for the root system $A_{m-1}$ can be determined as $I_{S_m/S_{m-1}}$.  The weight lattice of $A_{m-1}$, $\Lambda(A_{m-1})$ is its $\zz$-dual, and so can be identified with $J_{S_m/S_{m-1}}\cong \zz [S_m/S_{m-1}]/\zz$.

In particular, 
Klyachko determined $X_{mn}$ as a projective toric variety with character lattice  $M=\zz A_{m-1}\otimes \zz A_{n-1}$ given as the tensor product of root lattices for the root systems $A_{m-1}$ and $A_{n-1}$.  Hence its cocharacter lattice would be the tensor product of the corresponding weight lattices $N=\Lambda(A_{m-1})\otimes \Lambda(A_{n-1})$.

Let $\{e_i:1\le i\le m\}$ be the permutation basis for $\zz [S_m/S_{m-1}]$ and $\{\overline{e_i}:1\le i\le m\}$ its image in $J_{S_m/S_{m-1}}$ and $\{f_j:1\le j\le n\}$ be the permutation basis for $\zz [S_n/S_{n-1}]$ and $\{\overline{f_j}:j=1,\dots,n\}$ its image in $J_{S_n/S_{n-1}}$.  To simplify notation, we set 
$$\{e_{ij}=\overline{e_i}\otimes \overline{f_j} , 1\le i\le m, 1\le j\le n\}$$
and note that the above set is a  generating set for $J_{S_m/S_{m-1}}\otimes J_{S_{n}/S_{n-1}}$ (which can be identified with $N=\Lambda(A_{m-1})\otimes \Lambda(A_{n-1})$).
We observe that $\sum_{i=1}^me_{ij}=0$ and $\sum_{j=1}^ne_{ij}=0$ for all i,j.
One can choose $$\{e_{ij}:1\le i\le m-1, 1\le j\le n-1\}$$ as a $\zz$-basis of $N=\Lambda(A_{m-1})\otimes \Lambda(A_{n-1})$.

\begin{definition}
Let $$E_{mn}=\{e_{ij}:1\le i\le m, 1\le j\le n\}\subseteq N$$
We define $P_{mn}$ to  be the lattice polytope determined as the convex hull of $E_{mn}\subseteq N_{\rr}$.
Let $\Sigma_{mn}$ be the fan determined as cones over faces of $P_{mn}$. Then the Klyachko variety $X_{mn}$ is the toric variety corresponding to the fan $\Sigma_{mn}$. 
\end{definition}

The fan $\ss_{mn}$ comes from the faces of a polytope that has the origin as an interior point, therefore it is complete. The toric variety $X_{mn}$ is then complete, normal, and separated. Also, taking transposes of matrices gives a clear isomorphism $\ss_{mn}\simeq \ss_{nm}$ which implies that $X_{mn}\simeq X_{nm}$. Therefore, we can assume $n\geq m$ without loss of generality.

By this definition, the non-zero cones of $\ss_{mn}$ are in correspondence with proper faces of $P_{mn}$. This implies that we can explicitly write $\ss_{mn}$ as
$$\ss_{mn}=\{\cone(F) \:|\: F \text{ is a proper face of } P_{mn}\}\cup \{0\}$$
where $0$ is the cone consisting only of the origin.

The $e_{ij}$ are the vertices of $P_{mn}$ and are minimal elements of the lattice in the sense that there is no positive $c<1$ such that $ce_{ij}\in N$. This implies that the set of 1-dimensional cones is
$$\ss_{mn}(1)=\{\cone(e_{ij})\:|\: e_{ij}\in E_{mn}\}$$
and we have that $E_{mn}$ is the set of minimal generators $\{u_\rho\}_{\rho\in\ss_{mn}(1)}$ for the rays of $\ss_{mn}$. Specifically, we have $e_{ij}=u_{\cone(e_{ij})}$ and a bijection
\begin{align*}
E_{mn} &\longleftrightarrow \ss_{mn}(1) \\
e_{ij} &\longleftrightarrow \cone(e_{ij}).
\end{align*}
For this reason, we identify the elements of $\ss_{mn}(1)$ with the elements of $E_{mn}$ to which they correspond.

By construction, we note that $S_m\times S_n$ acts transitively on the vertices of $P_{mn}$ since $E_{mn}\subseteq N$ lie in one orbit of the natural action of $S_{m}\times S_n$ on $N=J_{S_{m}/S_{m-1}}\otimes J_{S_n/S_{n-1}}$. As any proper face of $P_{mn}$ is the convex hull of a subset of $E_{mn}$, we see that $S_m\times S_n$ also acts transitively on the faces of $P_{mn}$ and on the cones of the fan $\ss_{mn}$.

In particular, any cone of the fan $\Sigma_{mn}$, can be expressed as $\cone(I)$ for some $I\subset E_{mn}$. The action of $\tau\in S_{m}\times S_n$ on $E_{mn}\subseteq N$, then induces an action $\tau\cdot \cone(I)=\cone(\tau(I))$ on the cones of the fan $\Sigma_{mn}$.

The lattice polytope $P_{mn}$ has many other beautiful properties, which determine geometric properties of the Klyachko varieties. These were proven by Klyachko and Voskresenskii, and then in different, more explicit language by Codogni, Fanelli, Svaldi, and Tasin.

We first recall some polytope definitions from ~\cite{KN13}\footnote{Note that these definitions vary from source to source.}

\begin{definition}
A convex lattice polytope $P\subseteq N_{\rr}$, $\dim(P)=\rank(N)$ is called Fano if:
\begin{enumerate}
\item The origin is contained in the strict interior of $P$.
\item Each vertex of $P$ is a primitive lattice point of $P$.
\end{enumerate}
A Fano polytope $P$ is in particular smooth if the vertices of each facet form a basis for the lattice $N$.
\end{definition}

In particular, we have already observed that $P_{mn}$ is a Fano polytope.

\begin{definition}
A convex lattice polytope $P\subseteq N_{\rr}$, $\dim(P)=\rank(N)$ is called 

\begin{itemize} 
\item \emph{reflexive} if its polar dual polytope $$P^0=\{v\in M_{\rr}: \langle v,u\rangle \ge -1, u\in P\}$$
is also a lattice polytope.
\item \emph{terminal} if the only lattice points of $P$ are the origin and its vertices.
\end{itemize}

\end{definition}

A Gorenstein Fano toric variety $X_{\Sigma}$ is a complete normal toric variety whose anticanonical divisor is Cartier and ample. A Gorenstein Fano toric variety $X_{\Sigma}$ is  in particular projective, and is equipped with a projective embedding determined by a multiple of the anticanonical divisor $-K_{X_{\Sigma}}=\sum_{\rho \in \Sigma(1)}D_{\rho}$ determined by the rays of the fan.  

That a Gorenstein Fano toric variety $X_{\Sigma}$ has terminal singularities means roughly that it has a ``well-behaved'' resolution of singularities.  More precisely, if $f:X'\to X_{\Sigma}$ is a proper birational resolution of singularities of $X_{\Sigma}$, then one can describe the relationship between the  canonical divisor of $X'$, $K_{X'}$ and that of the pullback $f^*X_{\Sigma}$.  We would have 
$$K_{X'}=f^*K_{X_{\Sigma}}+\sum_{i=1}^ma_iE_i$$ where $E_i,i=1,\dots,m$ are the exceptional divisors for some $a_i\in \mathbb{Q}$.
If in addition, $X_{\Sigma}$ has terminal singularities, it is the case that $a_i\in \mathbb{Q}_{>0}$, for all $i$.

\begin{proposition}~\cite{VK84,CFST18}
$X_{mn}$ is a Gorenstein Fano toric variety which has terminal singularities.
\end{proposition}

A toric variety $X_{\Sigma}$ is Gorenstein Fano if and only if the polytope corresponding to its anticanonical divisor 
$$\{v\in M_{\rr}: \langle v,u_{\rho}\rangle \ge -1, \rho\in \Sigma(1)\}$$
is reflexive.
We remark that this is in turn equivalent to the fact that the Fano polytope generated by primitive lattice vectors $\conv\{u_{\rho}:\rho \in \Sigma(1)\}$ is reflexive, as these two polytopes are polar duals.

A Gorenstein Fano toric variety has terminal singularities if and only if its Fano polytope is terminal.

Note that the polar dual $P_{mn}$ may be identified with a translate of the central transportation polytope:

  $$\Omega_{m,n}=\{(y_{ij})\in M_{mn}(\rr):\sum_{i=1}^my_{ij}=m,\sum_{j=1}^ny_{ij}=n,y_{ij}\ge 0\}$$

We will focus on the case in which $X_{mn}$ is smooth. 

\begin{proposition}~\cite{VK84, CFST18}
The polytope $P_{mn}$ is a smooth Fano polytope if and only if $m$ and $n$ are relatively prime. In particular, this implies that the variety $X_{mn}$ is smooth if and only if $m$ and $n$ are relatively prime.
\end{proposition}

Some nice results involving the Klyachko varieties are given in\cite{VK84}. In particular, Corollary 2, Proposition, and the Corollary following Theorem 7 give us information about $P_{mn}$ and the rationality of forms of $X_{mn}$. 

Given a field extension $L/\fff$, an $L/\fff$ form of a variety $X$ is an $\fff$-variety $Y$ such that $X_L=Y_L$.  

\begin{proposition}\label{xmn smooth}
Suppose $m$ and $n$ are relatively prime.
Any form of $X_{mn}$ is stably rational and if $n= 1 \mod m$ then it is rational.
\end{proposition}

The Klyachko variety $X_{mn}$ can be used to determine a toric model of the (non-split) algebraic torus corresponding to the $S_m\times S_n$-lattice $\zz A_{m-1}\otimes \zz A_{n-1}$.\cite{Kly83,Lem}.

\subsection{Matrix Representations}

Let us now establish some notation.
Many of the properties of the variety $X_{mn}$ follow from properties of the fan $\ss_{mn}$. 

Let us fix the integers $m$ and $n$ so that we may write $X$, $\ss$, $P$, and $E$ instead of $X_{mn}$, $\ss_{mn}$, $P_{mn}$, and $E_{mn}$, respectively.

\begin{definition}
Let $I\subset E$ be a subset. We define a $m\times n$ matrix $[I_{ij}]$ where $I_{ij}=1$ if $e_{ij}\in I$ and $I_{ij}=0$ if $e_{ij}\not\in I$. We call the matrix $[I_{ij}]$ the {\it matrix representation} for $I$.
\end{definition}

For a subset $I\subset E$, we abuse notation and write $I$ for both the subset and its corresponding matrix representation. This correspondence gives a bijection
\begin{align*}
\ppp(E) &\longleftrightarrow M_{m,n}(\{0,1\}) \\
I\subset E &\longleftrightarrow I=[I_{ij}]
\end{align*}
and each entry $I_{ij}$ is the value of the indicator function
\begin{align*}
1_I:E&\to \{0,1\} \\
e_{ij}&\mapsto \begin{cases}
1 & \text{if } e_{ij}\in I \\
0 & \text{if } e_{ij}\notin I
\end{cases}
\end{align*}
at $e_{ij}$.

In particular, the matrix representation of the empty set, $\emptyset$ is the $m\times n$ zero matrix and the matrix representation of the set $E$ is the all ones $m\times n$ matrix.

This notation is convenient because the group action of $S_m\times S_n$ acts on the matrix representations by permuting the rows and columns. Specifically, if $\tau=\tau_1\times \tau_2\in S_m\times S_n$, we can compute $\tau I$ by permuting the rows of $I$ by $\tau_1$ and the columns by $\tau_2$.
This then determines the action on the cone corresponding to the subset $I$ of $E$ as
$$\tau\cdot \cone(I)=\cone(\tau\cdot I)$$
for all $\tau \in S_m\times S_n$ and $I\subset E$.

\section{Maximal Cones and Primitive Collections}

In this chapter, we will give a full description of the maximal cones and primitive collections for $\ss$ when $n=1 \mod m$. The definition of a primitive collection is given below.

\begin{definition}
A {\it primitive collection} for a fan $\ss$ is a minimal subset $P\subset\ss(1)$ such that $P$ is not contained in any cone of $\ss$.
\end{definition}

Note that this definition implies that if $P$ is a primitive collection for $\ss$, then any proper subset of $P$ is contained in some cone of $\ss$.

\section{Graph Theory}

The results of ~\cite{VK84} show that graph theory plays a key role in the fan $\ss_{mn}$. Their Proposition 5 and Corollary 1 (page 237) give a relation between the rays of $\ss_{mn}$ with the bipartite graph $\Ga_{mn}$ defined below.

\begin{definition}
Let $\Ga=\Ga_{mn}$ be a graph with vertices $e_1,\ldots,e_m,f_1,\ldots,f_n$ and edges $(e_i,f_j)$ for all $i,j$.
\end{definition}

The graph $\Ga$ is a complete bipartite graph with $n+m$ vertices and $mn$ edges. Moreover, the edges $(e_i,f_j)\in \Ga$ correspond to elements $e_{ij}\in E$. We let $V_e=\{e_1,\ldots,e_m\}$ and $V_f=\{f_1,\ldots,f_n\}$ be the disjoint sets that make up the vertices of $\Ga$ such that any edge in $\Ga$ connects a vertex in $V_e$ with one in $V_f$. We also let a set of vertices $V$ denote the graph $(V,\emptyset)$. In particular, $V_e$ and $V_f$ are used to denote the graphs $(V_e,\emptyset)$ and $(V_f,\emptyset)$, respectively.

In general, the order of a graph $G$, denoted $|G|$, is the cardinality of its set of vertices. If $T$ is a subset of edges of $G$, then $G\setminus T$ is the graph obtained from $G$ by deleting all of the edges in $T$. If $G_1=(V_1,E_1)$ and $G_2=(V_2,E_2)$ are graphs, then their union is $G_1\cup G_2=(V_1\cup V_2,E_1\cup E_2)$ and their intersection is $G_1\cap G_2=(V_1\cap V_2,E_1\cap E_2)$. The graphs $G_1$ and $G_2$ are considered disjoint if $V_1\cap V_2$ and $E_1\cap E_2$ are both empty. We say that $G_1$ is a subgraph of $G_2$ if $V_1\subset V_2$ and $E_1\subset E_2$ and we will write this as $G_1\subset G_2$.

\begin{definition}
Given a subset $I\subset E$, we construct a subgraph $\ib\subset\Ga$ as follows: $\ib$ will have the same set of vertices as $\Ga$ and we will have an edge $(e_i,f_j)\in \ib$ if and only if $e_{ij}\notin I$. We will call $\ib$ the complement graph of $I$.
\end{definition}

The edges in the complement graph $\ib$ come from elements of $E$ that are not in $I$. This gives a bijection between subsets of $E$ and the subgraphs of $\Ga$ with vertex set $V_e\cup V_f$
$$I\subset E \longleftrightarrow \ib \subset \Ga$$
where the bijection is defined by the relation
$$e_{ij}\in I^C \Longleftrightarrow (e_i,f_j)\in\ib.$$ In particular, we have that $\Ga$ is the complement graph of the empty set ($\overline{\emptyset}=\Ga$) and the complement graph of $E$ is the graph with vertices $V_e\cup V_f$ and no edges.

The following definition provides us with a useful function for relating the vertices $V_e$ and $V_f$ in the complement graphs.

\begin{definition}
For a subset $I\subset E$, we construct a function on subsets $S\subset V_e$. For such a set $S$, we let $f_I$ be the function
\begin{align*}
f_I:\ppp(V_e)&\to \ppp(V_f) \\
S &\mapsto \{f_j \:|\: (e_i,f_j)\in \ib \text{ for some } e_i\in S\}
\end{align*}
that assigns to $S$ the subset of $V_f$ of vertices that are directly connected to $S$.

We also define $d_I(S)$ to be the size of the set $f_I(S)$, ie. $d_I(S)=|f_I(S)|$.

Similarly, we define
\begin{align*}
e_I:\ppp(V_f)&\to \ppp(V_e) \\
S' &\mapsto \{e_i \:|\: (e_i,f_j)\in \ib \text{ for some } f_j\in S'\}
\end{align*}
and $d_I(S')=|e_I(S')|$ for subsets $S'\subset V_f$.
\end{definition}


To make notation more convenient, we let $e_i$ and $f_j$ also represent the sets $\{e_i\}$ and $\{f_j\}$, respectively. This way, for a fixed subset $I$, we may write $f_I(e_i)$ (resp. $e_I(f_j)$) instead of $f_I(\{e_i\})$ (resp. $e_I(\{f_j\})$) and we have that $d(e_i)=|f_I(e_i)|$ (resp. $d(f_j)=|e_I(f_j)|$) is the degree of the vertex $e_i$ (resp. $f_j$) in the complement graph $\ib$ since the complement graphs are bipartite.

\begin{proposition}\label{ficupcap}
For any sets $I\subset E$ and $S_1,S_2\subset V_e$, we have the relations
\begin{itemize}
\item[(1)] $f_I(S_1\cup S_2)=f_I(S_1)\cup f_I(S_2)$,
\item[(2)] $f_I(S_1\cap S_2)\subset f_I(S_1)\cap f_I(S_2)$.
\end{itemize}
\end{proposition}

\begin{proof}
We first prove (1). We have $f_j\in f_I(S_1\cup S_2)$ if and only if $(e_i,f_j)\in \ib$ for some $e_i\in S_1\cup S_2$. This is equivalent to the condition that $(e_i,f_j)\in \ib$ for some $e_i\in S_1$ or $e_i\in S_2$ which is the same as $f_j\in f_I(S_1)\cup f_I(S_2)$.

Now let $f_j\in f_I(S_1\cap S_2)$. Then there is some $e_i\in S_1\cap S_2$ with $(e_i,f_j)\in \ib$. Since $e_i\in S_1$ and $(e_i,f_j)\in \ib$, we have $f_j\in f_I(S_1)$. Similarly, we also have $f_j\in f_I(S_2)$ and therefore $f_j\in f_I(S_1)\cap f_I(S_2)$. This proves (2).
\end{proof}

Although there is an identical result for the functions $e_I$ on the subsets of $V_f$, this proposition is sufficient to prove our results. This is due to the fact that the set $V_f$ is generally much larger than $V_e$.

\section{A Sufficient Condition on Generating sets for Maximal Cones}

For the remainder of this paper, let $n=mk+1$ for integers $m,n,k$ with $m,n>1$, $k\geq 0$. The results of [VK85] imply that we need to look at spanning trees in $\Ga$ in order to find the maximal cones of $\ss$. Specifically, we know that $X$ is smooth by Proposition \ref{xmn smooth}. This means that generating sets $I\subset E$ for the maximal cones of $\ss$ must form a basis for $N$. Corollary 1 (page 237) of [VK85] tells us that bases of $N$ correspond to spanning trees of $\Ga$, however, it is still not clear which spanning trees correspond to generating sets of maximal cones.

Our next result helps remedy this, we give a criterion on the degree of the vertices of a spanning tree that ensures it is the complement graph of a generating set for some maximal cone in $\ss$.

\begin{proposition}\label{max cones}
For $n=mk+1$, any set $I$ such that $\ib$ is a tree and $d_I(e_i)=k+1$ for all $i$ will generate a maximal cone of $\ss_{mn}$.
\end{proposition}

\begin{proof}
Note that any tree in $\Ga_{mn}$ in which $d_I(e_i)=k+1$ for all $i$ has number of edges equal to
\begin{align*}
\sum_{i=1}^m d_I(e_i)&=m(k+1) \\
&=m+n-1
\end{align*}
which is the number of edges in a spanning tree of a graph with $m+n$ vertices. Therefore, these trees are always spanning trees.

Let $I\subset E$ be a subset such that $\ib$ is a tree with $d_I(e_i)=k+1$ for all $i$. By the above reasoning, $\ib$ is a spanning tree, so deleting any edge $(e_i,f_j)\in \ib$ gives a graph $\ib\setminus \{(e_i,f_j)\}$ that has two connected components. For each $(e_i,f_j)\in \ib$, we define $G_{ij}$ to be the connected component of $\ib\setminus \{(e_i,f_j)\}$ that contains the vertex $e_i$ and we let
$$g_{ij}=|G_{ij}\cap V_e|.$$

Let $j\in\lln$ be fixed. There must be at least one $e_i\in e_I(f_j)$ since $\ib$ is a spanning tree, so the set of graphs $\{G_{ij} \:|\: e_i\in e_I(f_j)\}$ is non-empty.

Since $\ib$ is a spanning tree, there is a unique path in $\ib$ from any vertex in $V_e$ to $f_j$ and this path must pass through one of the vertices in the set $e_I(f_j)$ since $\ib$ is bipartite. If the path passes through $e_i\in e_I(f_j)$, then $e_i$ is unique otherwise there would be a cycle. This implies that the original vertex (from which the path started) is a vertex in the corresponding graph $G_{ij}$ and thus
$$V_e\subset \cup_{e_i\in e_I(f_j)} G_{ij}.$$
Also, since the vertex $e_i$ which the path passes through is unique, we must have that the graphs in the set $\{G_{ij} \:|\: e_i\in e_I(f_j)\}$ are all disjoint.

One may also observe that the $G_{ij}$ for $e_i\in e_I(f_j)$ along with the vertex $f_j$ form the set of connected components of the graph $\ib\setminus \{(e_i,f_j) \:|\: e_i\in e_I(f_j)\}$. This observation makes the previous facts and the following implication easier to understand.

The above facts imply that
$$V_e=\sqcup_{e_i\in e_I(f_j)} (G_{ij}\cap V_e)$$
and taking orders gives us
$$m=\sum_{e_i\in e_I(f_j)} |G_{ij}\cap V_e|$$
which gives the relation
\begin{align}
\sum_{e_i\in e_I(f_j)} g_{ij} &=m. \label{gijm}
\end{align}

Now let $i\in\llm$ be fixed. Then $f_I(e_i)$ is non-empty since $\ib$ is a spanning tree so the set of graphs $\{G_{ij} \:|\: f_j\in f_I(e_i)\}$ is also non-empty.

Let $G_{ij}^C$ be the connected component of $\ib\setminus \{(e_i,f_j)\}$ containing $f_j$. The notation is due to the fact that it is the graph complement to $G_{ij}$ in $\ib\setminus \{(e_i,f_j)\}$ in the sense that $G_{ij}\cup G_{ij}^C=\ib\setminus \{(e_i,f_j)\}$.

We will use logic similar to what we used before, most of which should be clear if one sees that the graphs $G_{ij}^C$ for $f_j\in f_I(e_i)$ and the vertex $e_i$ are the connected components of the graph $\ib\setminus \{(e_i,f_j) \:|\: f_j\in f_I(e_i)\}$.

Since $\ib$ is a bipartite spanning tree, there is always a unique path from each vertex in $V_e\setminus\{e_i\}$ to $e_i$ passing through a unique $f_j\in f_I(e_i)$. Thus, each vertex of $V_e\setminus\{e_i\}$ is in exactly one of the $G_{ij}^C$. This implies that 
$$V_e\setminus\{e_i\}\subset \cup_{f_j\in f_I(e_i)} G_{ij}^C$$
and the graphs $G_{ij}^C\cap V_e$ for $f_j\in f_I(e_i)$ are disjoint. Since $e_i\in G_{ij}$ for all $G_{ij}$, we also have $e_i\notin G_{ij}^C$ for all $f_j\in f_I(e_i)$. Combining this information gives us
$$V_e\setminus\{e_i\}= \sqcup_{f_j\in f_I(e_i)} (G_{ij}^C\cap V_e)$$
and we can take the orders to get
\begin{align}
m-1&=\sum_{f_j\in f_I(e_i)} |G_{ij}^C\cap V_e|. \label{GijC}
\end{align}
By definitions, we have $\ib=G_{ij}\cup G_{ij}^C$ for all $f_j\in f_I(e_i)$. Since the union is disjoint, intersecting with $V_e$ and taking cardinalities gives
$$m=g_{ij}+|G_{ij}^C\cap V_e|.$$
We can use this in \eqref{GijC} to get
\begin{align*}
m-1&=\sum_{f_j\in f_I(e_i)} (m-g_{ij}) \\
&=m|f_I(e_i)|-\sum_{f_j\in f_I(e_i)} g_{ij}
\end{align*}
which we rearrange and use the relation $|f_I(e_i)|=d_I(e_i)=k+1$ to obtain
\begin{align}
\sum_{f_j\in f_I(e_i)} g_{ij}=n. \label{gijn}
\end{align}

Now we must prove that $I$ generates a maximal cone in $\ss$. We can do this by constructing a hyperplane $H$ that is a supporting hyperplane for the polytope $P_{mn}$ with $I=H\cap E$. We let
$$h_{ij}=\begin{cases}
1 & \text{if } e_{ij}\in I \\
1-g_{ij} & \text{if } e_{ij}\notin I \\
\end{cases}$$
for each $i\in\llm$, $j\in\lln$ and
$$h=[h_{ij}].$$

We will show that $h\in M$ by computing that the row and column sums are all $0$. This will give us a hyperplane
$$H=\{x\in \nr \:|\: \lb h, x\rb=1\}$$
in $\nr$ that satisfies the desired properties.

For a fixed $j$, we have
\begin{align*}
\sum_{i=1}^m h_{ij}&=\sum_{1\le i\le m, e_{ij}\in I}1+\sum_{1\le i\le m, e_{ij}\not\in I}(1-g_{ij}) \\
&=m- \sum_{e_i\in e_I(f_j)} g_{ij} \\
&=0
\end{align*}
by \eqref{gijm} which shows that the columns of $h$ sum to $0$. Then, fixing $i$ gives
\begin{align*}
\sum_{j=1}^n h_{ij}&= \sum_{1\le j\le n, e_{ij}\in I}1+\sum_{1\le j\le n, e_{ij}\not\in I}(1-g_{ij})\\
&=n- \sum_{f_j\in f_I(e_i)} g_{ij} \\
&=0
\end{align*}
by \eqref{gijn}, showing that the row sums are also $0$.

Now that we have shown $H$ is a hyperplane in $\nr$, it remains to show that $H$ that is a supporting hyperplane for $P_{mn}$ with $I=H\cap E$. It is clear from the definitions of $G_{ij}$ and $g_{ij}$ that we always have $g_{ij}>0$. This means we always have $h_{ij}=1$ for $e_{ij}\in I$ and $h_{ij}<1$ for $e_{ij}\not\in I$. The inner products are computed as $\langle h, e_{ij} \rangle=h_{ij}$ for each $e_{ij}\in E$, so the fact that $h_{ij}\leq 1$ for all $e_{ij}\in E$ implies that $H$ is a supporting hyperplane for $P_{mn}$. Since $h_{ij}=1$ exactly when $e_{ij}\in I$, we also have $I=H\cap E$.
\end{proof}

\section{The Maximal Cones and Primitive Collections}

It turns out that every maximal cone in $\ss$ comes from one of the sets described in Proposition \ref{max cones}. In order to prove this, we simultaneously consider the primitive collections for $\ss$. First, we show that the sets which we claim form primitive collections are not contained in any cone of $\ss$. Then if $I\subset E$, $I$ contains one of these sets or $I$ is part of a generating set for a maximal cone in $\ss$. We show that if the latter condition is satisfied, taking $I$ to be maximal gives one of the sets in Proposition \ref{max cones}.

We will need the following result.

\begin{lemma}\label{graphdeg}
A subset $I\subset \ss(1)$ does not contain any set $P$ of the form
$$P=\{e_{ij} \: |\: i\in Y,j\in Z\}$$
where $Y\subset \{1,\ldots m\}$ and $Z\subset \{1,\ldots ,n\}$ with $1\leq |Y|\leq m$ and $|Z|=(m-|Y|)k+1$ if and only if
\begin{align}
d_I(S)\geq |S|k+1
\end{align}
for any non-empty $S\subset V_e$.
\end{lemma}

\begin{proof}
Suppose $I\subset \ss(1)$ does not contain any of the sets $P$ described in the lemma, and let $A_i=\{j \:|\: e_{ij}\in I\}$ be the set of indices for the $i$-th row of $I$. Since $Z\subset \cap_{i\in Y} A_i$ if and only if $e_{ij}\in I$ for all $i\in Y$, $j\in Z$, the condition that $I$ does not contain $P=\{e_{ij} \: |\: i\in Y,j\in Z\}$ is equivalent to the condition
$$Z \not\subset \cap_{i\in Y} A_i.$$

For a given set $Y\subset \{1,\ldots m\}$ with $1\leq |Y|\leq m$, we must consider all $Z\subset \{1,\ldots ,n\}$ such that $|Z|=(m-|Y|)k+1$. Therefore, there cannot be a subset $Z\subset \cap_{i\in Y} A_i$ of size $(m-|Y|)k+1$, so we have
$$|\cap_{i\in Y} A_i| < (m-|Y|)k+1.$$
Since we must also consider all possible $Y\subset \{1,\ldots m\}$ with $1\leq |Y|\leq m$, the above inequality must hold for any such set $Y$. Taking the complement of the set on the left in this inequality then gives
$$|(\cap_{i\in Y} A_i)^C|> n-((m-|Y|)k+1)$$
which is reduced to
$$|\cup_{i\in Y} A_i^C| \geq |Y|k+1.$$

Now let $S\subset V_e$ be non-empty and consider $Y=\{i \:|\: e_i\in S\}$ so that $i\in Y$ if and only if $e_i\in S$ (hence $|Y|=|S|$). The previous inequality now implies that we have

\begin{align*}
d_I(S)&=|f_I(S)| \\
&=|\{f_j \:|\: (e_i,f_j)\in \ib \text{ for some } e_i\in S\}| \\
&=|\{j \:|\: e_{ij}\notin I \text{ for some } i\in Y\}| \\
&=|\cup_{i\in Y} A_i^C| \\
&\geq |Y|k+1 \\
&=|S|k+1.
\end{align*}
Since the other implication is obtained by reversing the argument, this completes the proof.
\end{proof}

We can now prove the main result of this chapter which is given in the following proposition.

\begin{proposition}\label{max cones2}
All of the maximal cones of $\ss$ are generated by subsets $I\subset E$ that correspond to trees in $\Ga$ in which $d_I(e_i)=k+1$ for all $i$.

The primitive collections for $\ss$ are the sets of the form
$$P=\{e_{ij} \: |\: i\in Y,j\in Z\}$$
where $Y\subset \{1,\ldots m\}$ and $Z\subset \{1,\ldots ,n\}$ with $1\leq |Y|\leq m$ and $|Z|=(m-|Y|)k+1$.
\end{proposition}

\begin{proof}
We begin by showing that these sets are not contained in any cone of $\ss$. Let $P=\{e_{ij} \: |\: i\in Y,j\in Z\}$ with $Y\subset \{1,\ldots m\}$, $Z\subset \{1,\ldots ,n\}$, $1\leq |Y|\leq m$, $|Z|=(m-|Y|)k+1$ and let $y=|Y|$ and $z=|Z|$. Suppose that this set is contained in a supporting hyperplane $H=\{x\in \nr \: |\: \lb h, x\rb =1\}$ where $h\in M$. This implies that we have

\begin{align*}
&h_{ij}\leq 1 & &\text{for all } i \text{ and } j \\
&h_{ij}=1 & &\text{for } i\in Y \text{ and } j\in Z \\
&\sum_{i=1}^m h_{ij}=0 & &\text{for all } j \\
&\sum_{j=1}^n h_{ij}=0 & &\text{for all } i.
\end{align*}

We set
$$A=\sum_{i\in Y^C}\sum_{j\in Z} h_{ij}$$
$$B=\sum_{i\in Y^C}\sum_{j\in Z^C} h_{ij}$$
and we have
\begin{align*}
A&=\sum_{i\in Y^C}\sum_{j\in Z} h_{ij} \\
&=\sum_{j\in Z}(\sum_{i=1}^{m} h_{ij}-\sum_{i\in Y} h_{ij}) \\
&=\sum_{j\in Z}\sum_{i=1}^{m} h_{ij}-\sum_{j\in Z}\sum_{i\in Y} h_{ij} \\
&=-yz.
\end{align*}
Similarly, we also have
\begin{align*}
B&=\sum_{i\in Y^C}\sum_{j\in Z^C} h_{ij} \\
&=\sum_{i\in Y^C}(\sum_{j=1}^{n} h_{ij}-\sum_{j\in Z} h_{ij}) \\
&=\sum_{i\in Y^C}\sum_{j=1}^{n} h_{ij}-\sum_{i\in Y^C}\sum_{j\in Z} h_{ij} \\
&=yz.
\end{align*}
However,
\begin{align*}
B&=\sum_{i\in Y^C}\sum_{j\in Z^C} h_{ij} \\
&\leq (m-y)(n-z) \\
&=(m-y)(yk)
\end{align*}
means we must have $yz\leq (m-y)(yk)$, but this is a contradiction since $$yz=y[(m-y)k+1]=(m-y)(yk)+y.$$

The second part of the proof is showing that any set $I$ not containing any of the sets in the proposition is contained in some cone in $\ss$. We will do this by showing that any maximal such set generates a maximal cone in $\ss$. Since all of the sets $P$ given in the proposition are minimal with respect to each other (any two distinct sets do not contain each other), this will not only imply that the sets in the proposition are primitive collections, but also that any primitive collection is of this form.

Let $I$ be a maximal set that does not contain any of the sets in the proposition and note that this also means that $I$ is a maximal set satisfying Lemma \ref{graphdeg}. To show that $I$ generates a maximal cone, we will prove that $\ib$ is a tree with $d_I(e_i)=k+1$ for all $i$ (Proposition \ref{max cones}).

To show that $d_I(e_i)=k+1$ for all $i$, let $e_i\in V_e$ be arbitrarily fixed and let $r=d_I(e_i)$. Then the set $\{j \:|\: f_j\in f_I(e_i)\}$ has $r$ elements which we label $j_1,\ldots,j_r$. We now have that $(e_i,f_j)\in \ib$ if and only if $j\in \{j_1,\ldots,j_r\}$ which also implies that $e_{ij}\notin I$ if and only if $j\in \{j_1,\ldots,j_r\}$. In particular, $e_{ij_l}\notin I$ for each $l=1\ldots, r$, so the sets $I_l=I\cup \{e_{ij_l}\}$ strictly contain $I$. We also have relations
\begin{align}
&f_{I_l}(e_i)=f_I(e_i)\setminus\{f_{j_l}\} \label{fil1} && \\
&f_{I_l}(S)=f_I(S) && \text{ for any subset } S\subset V_e \text{ not containing } e_i \label{fil2} \\
&d_{I_l}(S)\leq d_I(S)\leq d_{I_l}(S) +1 && \text{ for any subset } S\subset V_e \label{fil3}
\end{align}
for each $l=1,\ldots,r$, all of which arise from the fact that $\overline{I_l}=\ib\setminus \{(e_i,f_{j_l})\}$.

By the maximality of $I$, none of the sets $I_l$ satisfy Lemma \ref{graphdeg}. Therefore, for each $l=1,\ldots, r$, there must be some non-empty $S_l\subset V_e$ such that
\begin{align}
d_{I_l}(S_l)<|S_l|k+1. \label{fil4}
\end{align}
We must also have $e_i\in S_l$ for each $l$, otherwise \eqref{fil2} and Lemma \ref{graphdeg} give
\begin{align*}
d_{I_l}(S_l)=|f_{I_l}(S)|=|f_I(S)|=d_I(S_l)\geq |S_l|k+1
\end{align*}
which contradicts \eqref{fil4}. Furthermore, Lemma \ref{graphdeg}, \eqref{fil3}, and \eqref{fil4} imply that $d_I(S_l)=|S_l|k+1$ since
\begin{align*}
|S_l|k +1\leq d_I(S_l)\leq d_{I_l}(S_l) +1 <|S_l|k+2.
\end{align*}

We can now use induction to compute $d_I(\cap_{l\in \Psi} S_l)=|\cap_{l\in \Psi} S_l|k+1$ for subsets $\Psi\subset \{1,\ldots,r\}$. We have just seen that this is true when $|\Psi|=1$, so suppose that $\Psi$ has a size greater than one and that the induction hypothesis is true for all proper subsets of $\Psi$. That is, suppose $|\Psi|>1$ and that
\begin{align}
d_I(\cap_{l\in \Psi'} S_l)=|\cap_{l\in \Psi'} S_l|k+1 \label{fil5}
\end{align}
when $\Psi'\subsetneq \Psi$.

Now, let $\psi\in \Psi$ and set $\Psi^*=\Psi\setminus\{\psi\}$. We have
\begin{align*}
d_I(\cap_{l\in \Psi} S_l)&=|f_I(\cap_{l\in \Psi} S_l)| \\
&\leq |f_I(\cap_{l\in \Psi^*} S_l)\cap f_I(S_{\psi})|
\end{align*}
by Proposition \ref{ficupcap}. We can use the inclusion-exclusion formula to obtain that this is equal to
$$|f_I(\cap_{l\in \Psi^*} S_l)|+ |f_I(S_{\psi})|- |f_I(\cap_{l\in \Psi^*} S_l)\cup f_I(S_{\psi})|$$
and then use the induction hypothesis on the first two terms (which are equal to $d_I(\cap_{l\in \Psi^*} S_l)$ and $d_I(S_{\psi})$, respectively) and Proposition \ref{ficupcap} on the third term to rewrite this as
\begin{align*}
(|\cap_{l\in \Psi^*} S_l|k+1)+(|S_{\psi}|k+1)- |f_I((\cap_{l\in \Psi^*} S_l)\cup S_{\psi})|.
\end{align*}
We now have
\begin{align}
d_I(\cap_{l\in \Psi} S_l)\leq  (|\cap_{l\in \Psi^*} S_l|k+1)+(|S_{\psi}|k+1)- |f_I((\cap_{l\in \Psi^*} S_l)\cup S_{\psi})|. \label{fil6}
\end{align}
We can apply Lemma \ref{graphdeg} and again use the inclusion-exclusion formula as follows:
\begin{align*}
|f_I((\cap_{l\in \Psi^*} S_l)\cup S_{\psi})|&=d_I((\cap_{l\in \Psi^*} S_l)\cup S_{\psi}) \\
&\geq |(\cap_{l\in \Psi^*} S_l)\cup S_{\psi}|k+1 \\
&=(|\cap_{l\in \Psi^*} S_l|+|S_{\psi}|-|\cap_{l\in \Psi} S_l|)k+1.
\end{align*}
We can apply this to \eqref{fil6} to see that
\begin{align*}
d_I(\cap_{l\in \Psi} S_l)&\leq  (|\cap_{l\in \Psi^*} S_l|k+1)+(|S_{\psi}|k+1)- |f_I((\cap_{l\in \Psi^*} S_l)\cup S_{\psi})| \\
&\leq (|\cap_{l\in \Psi^*} S_l|k+1)+(|S_{\psi}|k+1)- [(|\cap_{l\in \Psi^*} S_l|+|S_{\psi}|-|\cap_{l\in \Psi} S_l|)k+1] \\
&=|(\cap_{l\in \Psi} S_l)|)k+1.
\end{align*}
We know that $\cap_{l\in \Psi} S_l$ is always non-empty because $e_i\in S_l$ for all $l$, so Lemma \ref{graphdeg} can be applied to get that $|(\cap_{l\in \Psi} S_l)|)k+1$ is also a lower bound for $d_I(\cap_{l\in \Psi} S_l)$, hence they must be equal.

We may now take $\Psi=\{1,\ldots,r\}$ and set $S'=\cap_{l=1}^r S_l$ to obtain the result
\begin{align}
d_I(S')=|S'|k+1. \label{fil7}
\end{align}
If $S'=e_i$, then we obtain the desired result that $d_I(e_i)=k+1$. We will show that this is the only possibility for $S'$.

Suppose that the set $S^*=S'\setminus e_i$ is non-empty. Then we have
\begin{align*}
|f_I(e_i)\cap f_I(S^*)|&=|f_I(e_i)|+ |f_I(S^*)|-|f_I(e_i)\cup f_I(S^*)| \\
&=|f_I(e_i)|+ |f_I(S^*)|-|f_I(e_i\cup S^*)| \\
&=|f_I(e_i)|+ |f_I(S^*)|-|f_I(S')| \\
&=d_I(e_i)+ d_I(S^*)-d_I(S')
\end{align*}
by the inclusion-exclusion formula and Proposition \ref{ficupcap}. Next, we apply Lemma \ref{graphdeg} and \eqref{fil7} to this relation to see that $f_I(e_i)\cap f_I(S^*)$ is not empty:
\begin{align*}
|f_I(e_i)\cap f_I(S^*)|&=d_I(e_i)+ d_I(S^*)-d_I(S') \\
&\geq (k+1)+(|S^*|k+1)-(|S'|k+1) \\
&=1.
\end{align*}
Therefore, there must be some $f_{j_l}\in f_I(e_i)\cap f_I(S^*)$.

Let $l$ be such that $f_{j_l}\in f_I(e_i)\cap f_I(S^*)$. The following sequence uses Proposition \ref{ficupcap}, \eqref{fil1}, \eqref{fil2}, and the fact that $f_{j_l}\in f_I(S^*)$:
\begin{align*}
f_{I_l}(S')&=f_{I_l}(e_i\cup S^*) \\
&=f_{I_l}(e_i)\cup f_{I_l}(S^*) \\
&=(f_{I}(e_i)\setminus\{f_{j_l}\}) \cup f_I(S^*) \\
&=f_{I}(e_i)\cup f_I(S^*) \\
&=f_I(S').
\end{align*}
Now we use this and \eqref{fil2} to get that, for any set $S$ containing $S'$, we have
\begin{align*}
f_{I_l}(S)&=f_{I_l}(S'\cup (S\setminus S')) \\
&=f_{I_l}(S')\cup f_{I_l}(S\setminus S') \\
&=f_{I}(S')\cup f_{I}(S\setminus S') \\
&=f_{I}(S'\cup (S\setminus S')) \\
&=f_{I}(S).
\end{align*}
However, taking $S=S_l$ gives a contradiction because this would imply
\begin{align*}
|S_l|k+1\leq d_{I}(S_l)=|f_{I}(S_l)|=|f_{I_l}(S_l)|=d_{I_l}(S_l)<|S_l|k+1
\end{align*}
by Lemma \ref{graphdeg} and \eqref{fil4}. Therefore, the set $S^*$ cannot be non-empty and we must have that $S'=e_i$ and $d_I(e_i)=k+1$ as mentioned earlier.

Now that we have proved $d_I(e_i)=k+1$ for all $i$, we must show that $\ib$ is a (spanning) tree. We begin this step by proving that it has only one connected component. Let $\{G_\phi\:|\: \phi\in \Phi\}$ be the set of connected components of the graph $\ib$ where $\Phi$ serves as an indexing set. Then for any distinct $\phi, \alpha\in \Phi$ we have
\begin{align*}
f_I(G_{\phi}\cap V_e)\cap f_I(G_{\alpha}\cap V_e)=\emptyset
\end{align*}
since $G_{\phi}$ and $G_{\alpha}$ are disjoint and this implies that the union $\cup_{\phi\in \Phi} f_I(G_{\phi}\cap V_e)$ is disjoint.

We also have that $V_e=\cup_{\phi\in \Phi} (G_{\phi}\cap V_e)$ since every $e_i$ is in some connected component of $\ib$. Since this union is also disjoint, we must have
$$\sum_{\phi\in \Phi} |G_{\phi}\cap V_e|=m.$$
Now we use these facts and apply Lemma \ref{graphdeg} and Proposition \ref{ficupcap} to get
\begin{align*}
n&=d_I(V_e) \\
&=|f_I(V_e)| \\
&=|f_I(\cup_{\phi\in\Phi} (G_{\phi}\cap V_e))| \\
&=|\cup_{\phi\in\Phi} f_I(G_{\phi}\cap V_e)) \\
&=\sum_{\phi\in\Phi} |f_I(G_{\phi}\cap V_e)| \\
&=\sum_{\phi\in\Phi} d_I(G_{\phi}\cap V_e) \\
&\geq \sum_{\phi\in\Phi} (|G_{\phi}\cap V_e|k+1) \\
&=(\sum_{\phi\in\Phi} |G_{\phi}\cap V_e|)k+\sum_{\phi\in\Phi} 1 \\
&=mk+|\Phi|
\end{align*}
which can only be true if $|\Phi|=1$.

Since $\ib$ is a bipartite graph, the number of edges in $\ib$ is equal to
\begin{align*}
\sum_{e_i\in V_e} d(e_i)&=m(k+1) \\
&=m+n-1
\end{align*}
which is one less than its number of vertices. Since we have already shown it is connected, this means it must be a tree.

\end{proof}

\begin{example}
It is known that when $m=2$, the maximal cones of $\ss_{2n}$ correspond to choosing two disjoint sets $A,B$ of $\{1,\ldots, n\}$ each of size $(n-1)/2$. In terms of the proposition, we can generate a graph satisfying $d(e_i)=k+1$ from such subsets by taking the edges
$$\{(e_1,f_j) \: |\: j\not\in A\}\cup \{(e_2,f_j) \: |\: j\not\in B\}.$$
In this case, the collection of maximal cones is transitive with respect to the action of $S_n$ and thus all maximal cones lie in a single $S_2\times S_n$-orbit.

Let us now look at the next smallest case: $m=3$ and $n=4$. We will use Proposition \ref{max cones2} to find a representative for each orbit of maximal cones and primitive collections for $\ss=\ss_{34}$. Note that $k=1$.

We begin by looking at spanning trees in the complete bipartite graph $\Ga=\Ga_{34}$ whose vertices are comprised of the sets $V_e=\{e_1,e_2,e_3\}$ and $V_f=\{f_1,f_2,f_3,f_4\}$. Let $T_1$ and $T_2$ be the trees in the following image.
\begin{center}
\begin{tikzpicture}
\path (1.5,6.5) node {$T_1$};
\draw[very thin] (2,0) -- (1,1) -- (2,2) -- (1,3) -- (2,4) -- (1,5) -- (2,6);
\filldraw[black] (2,0) circle (1pt) node[anchor=north]{$f_4$};
\filldraw[black] (1,1) circle (1pt) node[anchor=north]{$e_3$};
\filldraw[black] (2,2) circle (1pt) node[anchor=north]{$f_3$};
\filldraw[black] (1,3) circle (1pt) node[anchor=north]{$e_2$};
\filldraw[black] (2,4) circle (1pt) node[anchor=north]{$f_2$};
\filldraw[black] (1,5) circle (1pt) node[anchor=north]{$e_1$};
\filldraw[black] (2,6) circle (1pt) node[anchor=north]{$f_1$};

\path (6.5,6.5) node {$T_2$};
\draw[very thin] (7,0) -- (6,1) -- (7,6);
\draw[very thin] (7,2) -- (6,3) -- (7,6);
\draw[very thin] (7,4) -- (6,5) -- (7,6);
\filldraw[black] (7,0) circle (1pt) node[anchor=north]{$f_4$};
\filldraw[black] (6,1) circle (1pt) node[anchor=north]{$e_3$};
\filldraw[black] (7,2) circle (1pt) node[anchor=north]{$f_3$};
\filldraw[black] (6,3) circle (1pt) node[anchor=north]{$e_2$};
\filldraw[black] (7,4) circle (1pt) node[anchor=north]{$f_2$};
\filldraw[black] (6,5) circle (1pt) node[anchor=north]{$e_1$};
\filldraw[black] (7,6) circle (1pt) node[anchor=north]{$f_1$};
\end{tikzpicture}
\end{center}

The trees $T_1$ and $T_2$ are the only spanning trees in $\Ga$ that satisfy $d(e_i)=k+1=2$ for $i=1,2,3$ up to $S_3$ and $S_4$ permutations of the vertices $V_e$ and $V_f$, respectively. In other words, any other spanning tree of $\Ga$ with $d(e_i)=2$ for all $e_i$ can be obtained by relabeling the $e_i$ and/or $f_j$ in one of $T_1$ or $T_2$.

Let
\begin{align*}
A_1&=\{e_{ij} \: |\: (e_i,f_j)\not\in T_1\}, \\
A_2&=\{e_{ij} \: |\: (e_i,f_j)\not\in T_2\}
\end{align*}
be the sets corresponding to $T_1$ and $T_2$ under the bijection of subsets of $E$ with subgraphs of $\Ga$ coming from taking complement graphs. Under this bijection, we have $T_1=\overline{A_1}$ and $T_2=\overline{A_2}$ and the matrix representations of $A_1$ and $A_2$ look like
\begin{align*}
&
A_1=\begin{bmatrix}
0&0&1&1 \\
1&0&0&1 \\
1&1&0&0
\end{bmatrix},
&&
A_2=\begin{bmatrix}
0&0&1&1 \\
0&1&0&1 \\
0&1&1&0
\end{bmatrix}.
\end{align*}
By Proposition \ref{max cones2}, the sets $A_1$ and $A_2$ are generating sets for maximal cones in $\ss$, and the uniqueness of the trees $T_1$ and $T_2$ with respect to the action of $S_3\times S_4$ implies that $\cone(A_1)$ and $\cone(A_2)$ are representatives for the orbits of maximal cones.

Proposition \ref{max cones2} also tells us that the primitive collections are of the form $P=\{e_{ij} \: |\: i\in Y,j\in Z\}$ where $Y\subset \{1,\ldots 3\}$ and $Z\subset \{1,\ldots ,4\}$ with $1\leq |Y|\leq 3$ and $|Z|=4-|Y|$. The action of $S_3\times S_4$ will act on such a set simply by changing the indexing sets, but not their cardinalities. This implies that we can find a representative for each orbit of primitive collections by taking $Y=\{1,\ldots, y\}$ and $Z=\{1,\ldots, 4-y\}$ for each possible value of $y=1,2,3$. This gives the three representatives
\begin{align*}
P_1&=\{e_{ij}\: |\: i\in \{1\}, j\in \{1,2,3\}\}=\begin{bmatrix}
1&1&1&0 \\
0&0&0&0 \\
0&0&0&0
\end{bmatrix}, \\
P_2&=\{e_{ij}\: |\: i\in \{1,2\}, j\in \{1,2\}\}=\begin{bmatrix}
1&1&0&0 \\
1&1&0&0 \\
0&0&0&0
\end{bmatrix}, \\
P_3&=\{e_{ij}\: |\: i\in \{1,2,3\}, j\in \{1\}\}=\begin{bmatrix}
1&0&0&0 \\
1&0&0&0 \\
1&0&0&0
\end{bmatrix}
\end{align*}
of orbits of primitive collections.
\end{example}

In fact, we can extend the information that we have found in this example to the more general case that $m=3$ and $n=3k+1$. This is because any bipartite spanning tree on the vertex sets $V_e=\{e_1,e_2,e_3\}$ and $V_f=\{f_1,\ldots,f_n\}$ with $d(e_i)=k+1$ for all $i$ can be obtained by extending (a relabeling of) one of $T_1,T_2$ by equivariantly connecting the vertices $f_5,\ldots, f_n$ to the $e_i$. The trees obtained from extending $T_1$ and $T_2$ in this way look like the trees $T_1^*$ and $T_2^*$, respectively (up to a relabeling of the vertices $f_5,\ldots, f_n$), given below.

\begin{center}
\begin{tikzpicture}
\path (1.5,6.5) node {$T_1^*$};
\draw[very thin] (2,0) -- (1,1) -- (2,2) -- (1,3) -- (2,4) -- (1,5) -- (2,6);
\filldraw[black] (2,0) circle (1pt) node[anchor=west]{$f_4$};
\filldraw[black] (1,1) circle (1pt) node[anchor=east]{$e_3$};
\filldraw[black] (2,2) circle (1pt) node[anchor=west]{$f_3$};
\filldraw[black] (1,3) circle (1pt) node[anchor=east]{$e_2$};
\filldraw[black] (2,4) circle (1pt) node[anchor=west]{$f_2$};
\filldraw[black] (1,5) circle (1pt) node[anchor=east]{$e_1$};
\filldraw[black] (2,6) circle (1pt) node[anchor=west]{$f_1$};

\draw[very thin] (1,5) -- (3,5.5);
\filldraw[black] (3,5.5) circle (1pt) node[anchor=west]{$f_5$};
\draw[very thin] (1,5) -- (3,4.5);
\filldraw[black] (3,4.5) circle (1pt) node[anchor=west]{$f_{4+(k-1)}$};
\filldraw[black] (3,5.125) circle (0.3pt);
\filldraw[black] (3,5) circle (0.3pt);
\filldraw[black] (3,4.875) circle (0.3pt);

\draw[very thin] (1,3) -- (3,3.5);
\filldraw[black] (3,3.5) circle (1pt) node[anchor=west]{$f_{5+(k-1)}$};
\draw[very thin] (1,3) -- (3,2.5);
\filldraw[black] (3,2.5) circle (1pt) node[anchor=west]{$f_{4+2(k-1)}$};
\filldraw[black] (3,3.125) circle (0.3pt);
\filldraw[black] (3,3) circle (0.3pt);
\filldraw[black] (3,2.875) circle (0.3pt);

\draw[very thin] (1,1) -- (3,1.5);
\filldraw[black] (3,1.5) circle (1pt) node[anchor=west]{$f_{5+2(k-1)}$};
\draw[very thin] (1,1) -- (3,0.5);
\filldraw[black] (3,0.5) circle (1pt) node[anchor=west]{$f_n$};
\filldraw[black] (3,1.125) circle (0.3pt);
\filldraw[black] (3,1) circle (0.3pt);
\filldraw[black] (3,0.875) circle (0.3pt);


\path (6.5,6.5) node {$T_2^*$};
\draw[very thin] (7,0) -- (6,1) -- (7,6);
\draw[very thin] (7,2) -- (6,3) -- (7,6);
\draw[very thin] (7,4) -- (6,5) -- (7,6);
\filldraw[black] (7,0) circle (1pt) node[anchor=west]{$f_4$};
\filldraw[black] (6,1) circle (1pt) node[anchor=east]{$e_3$};
\filldraw[black] (7,2) circle (1pt) node[anchor=west]{$f_3$};
\filldraw[black] (6,3) circle (1pt) node[anchor=east]{$e_2$};
\filldraw[black] (7,4) circle (1pt) node[anchor=west]{$f_2$};
\filldraw[black] (6,5) circle (1pt) node[anchor=east]{$e_1$};
\filldraw[black] (7,6) circle (1pt) node[anchor=west]{$f_1$};

\draw[very thin] (6,5) -- (8,5.5);
\filldraw[black] (8,5.5) circle (1pt) node[anchor=west]{$f_5$};
\draw[very thin] (6,5) -- (8,4.5);
\filldraw[black] (8,4.5) circle (1pt) node[anchor=west]{$f_{4+(k-1)}$};
\filldraw[black] (8,5.125) circle (0.3pt);
\filldraw[black] (8,5) circle (0.3pt);
\filldraw[black] (8,4.875) circle (0.3pt);

\draw[very thin] (6,3) -- (8,3.5);
\filldraw[black] (8,3.5) circle (1pt) node[anchor=west]{$f_{5+(k-1)}$};
\draw[very thin] (6,3) -- (8,2.5);
\filldraw[black] (8,2.5) circle (1pt) node[anchor=west]{$f_{4+2(k-1)}$};
\filldraw[black] (8,3.125) circle (0.3pt);
\filldraw[black] (8,3) circle (0.3pt);
\filldraw[black] (8,2.875) circle (0.3pt);

\draw[very thin] (6,1) -- (8,1.5);
\filldraw[black] (8,1.5) circle (1pt) node[anchor=west]{$f_{5+2(k-1)}$};
\draw[very thin] (6,1) -- (8,0.5);
\filldraw[black] (8,0.5) circle (1pt) node[anchor=west]{$f_n$};
\filldraw[black] (8,1.125) circle (0.3pt);
\filldraw[black] (8,1) circle (0.3pt);
\filldraw[black] (8,0.875) circle (0.3pt);
\end{tikzpicture}
\end{center}

To prove this claim, we will use the following lemma.

\begin{lemma}\label{df1}
Let $n=mk+1$ and suppose $T$ is a spanning tree for $\Ga_{mn}$. Then the set of vertices $\{f_j \: |\: d(e_i)>1\}$ has size at most $m-1$.
\end{lemma}

\begin{proof}

Since $T$ is a spanning tree for $\Ga_{mn}$, we know that $T$ has $m+n-1$ edges, each connected to exactly one of the $f_j$, and $d(f_j)\geq 1$ for all $j$. These facts give
\begin{align*}
m+n-1&=\sum_j d(f_j) \\
&=\sum_{d(f_j)>1} d(f_j) + \sum_{d(f_j)=1} d(f_j) \\
&\geq 2|\{d(f_j)>1\}|+|\{d(f_j)=1\}| \\
&=n+|\{d(f_j)>1\}|
\end{align*}
which implies that the set $\{f_j \: |\: d(f_j)>1\}$ has size at most $m-1$.
\end{proof}

Let $T\subset \Ga_{3n}$ be a tree with $d(e_i)=k+1$ for all $i$ and we will show that $T$ is a relabeling of one of $T_1^*$ or $T_2^*$. By Lemma \ref{df1}, the set of vertices
$$\{f_j \: |\: d(f_j)>1\}$$
consists of either one or two elements.

First, suppose that $d(f_j)>1$ for only one value of $j$. In this case, we assume that $d(f_1)>1$ and $d(f_j)=1$ for $j\neq 1$ and we must have $d(f_1)=3$. This implies that $f_1$ is connected to each of $e_1,e_2,e_3$. Since the other $f_j$ are connected to one $e_i$ in such a way that $d(e_i)=k+1$ for each $i$, we must have that $T$ is (a relabeling of) $T_2^*$.

If there are exactly two $f_j$ with $d(f_j)>1$, then we can assume that they are $f_2$ and $f_3$ without loss of generality (we will construct $T_1^*$ in this case). The only possible degrees for these vertices are $d(f_2)=d(f_3)=2$ and we can assume that the vertices connected to $f_2$ are $e_1$ and $e_2$. Now, it cannot be the case that $f_3$ is connected to both $e_1$ and $e_2$, otherwise there exists a cycle in $T$. This implies that $f_3$ is connected to $e_3$ and exactly one of $e_1,e_2$. Again, without loss of generality, we assume that $f_3$ is connected to $e_2$. Since all other $f_j$ have $d(f_j)=1$, the only way we can have $d(e_i)=k+1$ for all $i$ is if $T$ looks like $T_1^*$ (up to a relabeling the $f_j$).

The subsets $A_1^*$ and $A_2^*$ of $E$ corresponding to $T_1^*$ and $T_2^*$ are the sets
\begin{align*}
A_1^*&=\{e_{ij} \: |\: (e_i,f_j)\not\in T_1^*\}, \\
A_2^*&=\{e_{ij} \: |\: (e_i,f_j)\not\in T_2^*\}
\end{align*}
and the matrix representations of $A_1^*$ and $A_2^*$ look like
\begin{align*}
A_1^*&=\left[\begin{array}{c|ccc}
A_1 & \myunderbrace{\begin{array}{ccc}
0&\cdots & 0 \\
1&\cdots & 1 \\
1&\cdots & 1
\end{array}}{k-1} & \myunderbrace{\begin{array}{ccc}
1&\cdots & 1 \\
0&\cdots & 0 \\
1&\cdots & 1
\end{array}}{k-1} & \myunderbrace{\begin{array}{ccc}
1&\cdots & 1 \\
1&\cdots & 1 \\
0&\cdots & 0
\end{array}}{k-1}
\end{array}\right], \vphantom{\left[ \underbrace{\begin{array}{ccc}
0&\cdots & 0 \\
1&\cdots & 1 \\
1&\cdots & 1
\end{array}}_{k-1}\right]} \\
A_2^*&=\left[\begin{array}{c|ccc}
A_2 & \myunderbrace{\begin{array}{ccc}
0&\cdots & 0 \\
1&\cdots & 1 \\
1&\cdots & 1
\end{array}}{k-1} & \myunderbrace{\begin{array}{ccc}
1&\cdots & 1 \\
0&\cdots & 0 \\
1&\cdots & 1
\end{array}}{k-1} & \myunderbrace{\begin{array}{ccc}
1&\cdots & 1 \\
1&\cdots & 1 \\
0&\cdots & 0
\end{array}}{k-1}
\end{array}\right]. \vphantom{\left[ \underbrace{\begin{array}{ccc}
0&\cdots & 0 \\
1&\cdots & 1 \\
1&\cdots & 1
\end{array}}_{k-1}\right]}
\end{align*}

Now Proposition \ref{max cones2} tells us that the sets $A_1^*$ and $A_2^*$ are generating sets for maximal cones in $\ss$, and the uniqueness of the trees $T_1^*$ and $T_2^*$ with respect to the action of $S_3\times S_n$ implies that $\cone(A_1^*)$ and $\cone(A_2^*)$ are representatives for the orbits of maximal cones.

Note that the last $3(k-1)$ columns of $A_1^*$ and $A_2^*$ can be rearranged to form $k-1$ copies of $I^C=\begin{bmatrix}
0&1&1 \\
1&0&1 \\
1&1&0
\end{bmatrix}$. This implies that we can choose $\tau\in S_3\times S_n$ corresponding to a permutation of columns so that
\begin{align*}
\tau A_1^*&=\left[\begin{array}{c|c|c|c}
A_1 & I^C & \cdots & I^C
\end{array}\right], \\
\tau A_2^*&=\left[\begin{array}{c|c|c|c}
A_2 & I^C & \cdots & I^C
\end{array}\right]
\end{align*}
are also generating sets for maximal cones, and we could alternatively take these cones as representatives for the orbits (where $I^C$ appears in each matrix $k-1$ times).

Finding representatives of the orbits of the primitive collections of $\ss_{3n}$ can be done in a similar way to the case when $n=4$. We consider the sets $Y=\{1,\ldots ,y\}$ and $Z=\{1,\ldots ,(3-y)k+1\}$ for each value of $y=1,2,3$ to obtain the matrix representations

\begin{align*}
P_1^*&=\left[ \begin{array}{cc}
\myunderbrace{\begin{matrix}
1&\cdots & 1 \\
0&\cdots & 0 \\
0&\cdots & 0\end{matrix}}{2k+1}&
\myunderbrace{\begin{matrix}
0&\cdots & 0 \\
0&\cdots & 0 \\
0&\cdots & 0
\end{matrix}}{k}
\end{array} \right], \vphantom{\left[\underbrace{\begin{matrix}
1&\cdots & 1 \\
0&\cdots & 0 \\
0&\cdots & 0\end{matrix}}_{2k+1}\right]} \\
P_2^*&=\left[ \begin{array}{cc}
\myunderbrace{\begin{matrix}
1&\cdots & 1 \\
1&\cdots & 1 \\
0&\cdots & 0\end{matrix}}{k+1}&
\myunderbrace{\begin{matrix}
0&\cdots & 0 \\
0&\cdots & 0 \\
0&\cdots & 0
\end{matrix}}{2k}
\end{array} \right], \vphantom{\left[\underbrace{\begin{matrix}
1&\cdots & 1 \\
0&\cdots & 0 \\
0&\cdots & 0\end{matrix}}_{2k+1}\right]} \\
P_3^*&=\left[ \begin{array}{cc}
\begin{matrix}
1 \\
1 \\
1
\end{matrix}&
\myunderbrace{\begin{matrix}
0&\cdots & 0 \\
0&\cdots & 0 \\
0&\cdots & 0
\end{matrix}}{3k}
\end{array}\right]. \vphantom{\left[\underbrace{\begin{matrix}
1&\cdots & 1 \\
0&\cdots & 0 \\
0&\cdots & 0\end{matrix}}_{2k+1}\right]}
\end{align*}

This example shows that, in the case that $m=3$ and $n=3k+1$, the orbits of maximal cones and primitive collections for the fan $\ss_{3n}$ are in bijection with the orbits of maximal cones and primitive collections for the fan $\ss_{34}$. To state this specifically, we let $O_M(\ss_{mn})$ be the set of $S_m\times S_n$-orbits of maximal cones in the fan $\ss_{mn}$ and $O_P(\ss_{mn})$ the set of $S_m\times S_n$-orbits of primitive collections of $\ss_{mn}$. We have seen that there are bijections
\begin{align*}
O_M(\ss_{3n})\simeq O_M(\ss_{34}), \\
O_P(\ss_{3n})\simeq O_P(\ss_{34})
\end{align*}
coming from the correspondence of representatives given by $\cone(A_i^*)\leftrightarrow \cone(A_i)$ and $P_i^*\leftrightarrow P_i$.

After seeing this correspondence, the question arises of whether there are bijections
\begin{align*}
O_M(\ss_{mn})\simeq O_M(\ss_{m,m+1}), \\
O_P(\ss_{mn}))\simeq O_P(\ss_{m,m+1})
\end{align*}
for all $m$ and $n=1 \mod m$.

\section{Number of \texorpdfstring{$S_m \times S_n$-Orbits} {TEXT}}

It is easy to see that $O_P(\ss_{mn}))\simeq O_P(\ss_{m,m+1})$ is always guaranteed by the fact that representatives for the orbits are in correspondence with sets $Y=\{1,\ldots,y\}$ for each $y=1,\ldots m$ (and hence there are always $m$ $S_m \times S_n$-orbits of primitive collections). Unfortunately, the bijection $O_M(\ss_{mn})\simeq O_M(\ss_{m,m+1})$ only holds for $m=2$ and $m=3$. However, there is still a generalization of the bijection $O_M(\ss_{3n})\simeq O_M(\ss_{34})$ for $m\geq 3$ which is given in the next proposition. In fact, if one first finds (a generating set for) a representative of each orbit of maximal cones in the fan $\ss_{m,m(m-2)+1}$, it is relatively easy to use this information to find a representative for each orbit of maximal cones and primitive collections for $\ss_{mn}$ anytime $n=mk+1$. This result will make use of the following lemma.

\begin{lemma}\label{zvf}
Let $m\ge 3$.
Suppose $n=mk+1$ with $a=k-(m-2)>0$ and that $A^*$ is a generating set for a maximal cone $\s\in \ss_{mn}$. Then there is some set $Z\subset V_f$ of size $ma$ such that $d_{A^*}(f_j)=1$ for all $f_j\in Z$ and $|Z\cap f_{A^*}(e_i)|=a$ for all $i$.
\end{lemma}

\begin{proof}
Let $Q=\{f_j \: |\: d_{A^*}(f_j)=1\}$ and suppose that $|f_{A^*}(e_{i_1})\cap Q|<a$ for some $i_1$. Then, since the sets $f_{A^*}(e_i)\cap Q$ are disjoint and their union is $Q$, we have
\begin{align*}
|Q|&=\sum_{i=1}^m |f_{A^*}(e_i)\cap Q| \\
&<a+\sum_{i\neq i_1} |f_{A^*}(e_i)\cap Q|.
\end{align*}
We also have
\begin{align*}
|Q|&=|V_f|- |\{f_j \: |\: d_{A^*}(e_i)>1\}| \\
&\geq n-(m-1)
\end{align*}
by Lemma \ref{df1} and the fact that $Q=V_f\setminus \{f_j \: |\: d_{A^*}(e_i)>1\}$ and we can combine these inequalities to obtain
$$n-(m-1)<a+\sum_{i\neq i_1} |f_{A^*}(e_i)\cap Q|.$$
We now rearrange and simplify this inequality as follows:
\begin{align*}
\sum_{i\neq i_1} |f_{A^*}(e_i)\cap Q|&>n-(m-1)-a \\
&=n-(m-1)-[k-(m-2)] \\
&=n-k-1 \\
&=(m-1)k.
\end{align*}
Since $\frac{(m-1)k}{m-2}>k$, at least one of the sets $f_{A^*}(e_i)\cap Q$ with $i\neq i'$ must have size at least $k+1$. This implies that there is some $i_2\neq i_1$ such that
\begin{align*}
k+1&\leq |f_{A^*}(e_{i_2})\cap Q| \\
&\leq |f_{A^*}(e_{i_2})| \\
&= d_{A^*}(e_{i_2}) \\
&=k+1
\end{align*}
where the last inequality follows from Proposition \ref{max cones2}. This can only be true if $f_{A^*}(e_{i_2})\cap Q=f_{A^*}(e_{i_2})$ which means that $f_{A^*}(e_{i_2})\subset Q$. However, this would imply that all vertices connected to $e_{i_2}$ have degree one and there would be no paths from $e_{i_2}$ to any other $e_i$, contradicting that $\overline{A^*}$ is a spanning tree. Thus, our original assumption that there is some $i_1$ for which $|f_{A^*}(e_{i_1})\cap Q|<a$ is impossible and we must have
$$|f_{A^*}(e_{i})\cap Q|\geq a \quad \text{for all } i.$$
This means that there exist disjoint subsets $Z_i\subset f_{A^*}(e_{i})\cap Q$ of size $a$ and we can take $Z=\cup_{i=1}^m Z_i$ to achieve the desired result.
\end{proof}

\begin{proposition}
If $n=mk+1$ with $m\geq 3$ and $k\geq m-2$, there is a bijection
\begin{align*}
O_M(\ss_{mn})\simeq O_M(\ss_{m,m(m-2)+1})
\end{align*}
described as follows. If $O(\s)\in O_M(\ss_{m,m(m-2)+1})$ is the orbit of a maximal cone $\s\in \ss_{m,m(m-2)+1}$ generated by a set $A\subset E_{m,m(m-2)+1}$ (ie., $\s=\cone(A)$), then the corresponding orbit $O(\s^*)\in O_M(\ss_{mn})$ can be found by taking $\s^*=\cone(A^*)$ where
$$A^*=\left[\begin{array}{c|c|c|c}
A&I^C &\cdots &I^C
\end{array}\right]$$
is constructed from the matrix representation for $A$ by adding $k-(m-2)$ copies of $I^C$.

\end{proposition}

\begin{proof}
Let $n'=m(m-2)+1$ and $a=k-(m-2)$. If $k=m-2$, the result is trivial, so suppose $a>0$.

First, we show that given $O(\s)\in O_M(\ss_{m,m(m-2)+1})$, the construction in the proposition gives a corresponding orbit $O(\s^*)\in O_M(\ss_{mn})$. We leave the proof that this induces a well-defined and one-to-one map $O_M(\ss_{m,m(m-2)+1})\to O_M(\ss_{mn})$ to the reader. Since the complement graph $\overline{A^*}$ is obtained from $\overline{A}$ by adding the $n-n'$ vertices $f_{n'+1},\ldots,f_n$ and, for each $j=n'+1,\ldots, n$, we add the edge $(e_{i_j},f_j)$ where $i_j$ is the residue of $j-n'$ modulo $m$, ie., $i_j=j-n' \mod m$ ($=j-1 \mod m$). Then we have
\begin{align*}
d_{A^*}(e_i)&=d_{A}(e_i)+(n-n')/m \\
&=d_{A}(e_i)+a \\
&=(m-1)+a \\
&=k+1
\end{align*}
which arises from the facts that there are $n-n'$ new edges equivariantly connected to the $e_i$ (resulting in (n-n')/m new edges connected to each $e_i$) and that $d_{A}(e_i)=(m-2)+1$ by Proposition \ref{max cones2}. Furthermore, it is clear that $\overline{A^*}$ is also a tree since each new edge is formed from connecting a new vertex to the graph $\overline{A}$.

Now, we show that every orbit of maximal cones $O(\s^*)\in O_M(\ss_{mn})$ comes from some orbit $O(\s)\in O_M(\ss_{m,m(m-2)+1})$ under the given construction. Suppose that $A^*$ is a generating set for some maximal cone in $\ss_{mn}$.

By Proposition \ref{max cones2}, the complement graph $\overline{A^*}$ is a tree with $d_{A^*}(e_i)=k+1$ and, as we have seen previously, this is a spanning tree for $\Ga_{mn}$ since the total number of edges is
\begin{align*}
\sum_j d_{A^*} (f_j)&=\sum_i d_{A^*} (e_i) \\
&=m(k+1) \\
&=m+n-1.
\end{align*}
By Lemma \ref{df1}, the set $\{j \: |\: d_{A^*} (f_j)>1\}$ has size at most $m-1$. Without loss of generality, we can assume
$$\{j \: |\: d_{A^*} (f_j)>1\}=\{1,\ldots ,r\}$$
for some $r\leq m-1$.

By Lemma \ref{zvf}, there is some set $Z\subset V_f$ of size $n-n'$ such that the set $d_{A^*}(f_j)=1$ for all $f_j\in Z$ and $|Z\cap f_{A^*}(e_i)|=a$ for all $i$. Then, without loss of generality, we can assume that $Z=f_{n'+1},\ldots,f_n$ and that
$$Z\cap f_{A^*}(e_i)=\{f_j\:|\: n'+1\leq j \leq n \text{ and } j=n'+i \mod m\}.$$ Letting $A\subset E$ be the set corresponding to the graph obtained from $\overline{A^*}$ by removing the vertices in $Z$ and the $n-n'$ edges connected to them will then give a generating set for a maximal cone of $\ss_{mn'}$. We prove this by showing that this new graph $\overline{A}\subset \Ga_{mn'}$ satisfies the conditions given in Proposition \ref{max cones2}.

Since all of the vertices in $Z$ are endpoints of $\overline{A^*}$, removing them and their edges preserves connectedness. Also, any cycle in $\overline{A}$ would give a cycle in $\overline{A^*}$, so there must not exist any since $\overline{A^*}$ is a tree. Therefore, $\overline{A}$ is a tree and the degrees of the $e_i$ vertices can be computed as
\begin{align*}
d_{A}(e_i)&=d_{A^*}(e_i)-|Z\cap f_{A^*}(e_i)| \\
&=k+1-a \\
&=(m-2)+1.
\end{align*}
By Proposition \ref{max cones2}, we have that $A$ is a generating set for a maximal cone of $\ss_{mn'}$ and $A^*$ is constructed from $A$ as described in the proposition.

\end{proof}

In the special case where $n=m+1$, generating sets for maximal cones of $\ss_{m,m+1}$ correspond to bipartite spanning trees for which $d(e_i)=2$ for all $i$. In this case, we can think of the two edges connected to each $e_i$ as a single edge connecting two of the $f_j$ vertices together. This idea leads to the following result.

\begin{proposition}
Let $U(n)$ be the set of unlabeled spanning trees on $n$ vertices.

The $S_m$-orbits of maximal cones of $\ss_{m,m+1}$ correspond to labeled spanning trees on $m+1$ vertices, and the $S_m\times S_{m+1}$-orbits correspond to unlabeled trees on $m+1$ vertices. Consequently,
$$|O_M(\ss_{m,m+1})|=|U(m+1)|$$
where $U(n)$ is the set of unlabeled spanning trees on $n$ vertices.

The number of unlabeled spanning trees on $n$ vertices, $|U(n)|$ (Sequence A000051 in \cite{OEIS}), can be determined from the generating function
$$A(x)=1+g(x)-g^2(x)/2+g(x^2)/2$$
where $g(x)=x+x^2+2x^3+\ldots$ is the generating function for the number of unlabeled rooted trees with n nodes (Sequence A000081 in \cite{OEIS}).
\end{proposition}

\begin{proof}
By Proposition \ref{max cones2}, we know that maximal cones of $\ss_{m,m+1}$ are in bijection with the set $T_{m,m+1;2}$ of bipartite trees on vertex sets $V_e$ of size $m$ and $V_f$ of size $m+1$ such that all vertices in $V_e$ have degree 2 and that this bijection is equivariant with respect to the $S_m\times S_{m+1}$ action. So, we have  $|O_M(\ss_{m,m+1})|=|T_{m,m+1;2}/(S_m\times S_{m+1})|$.

Let $T_{m+1}$ be the set of labeled trees on $m+1$ vertices.
We define a map $\Phi:T_{m,m+1;2}\to T_{m+1}$.  Each bipartite tree $T$ in $T_{m,m+1;2}$ is determined by its edge set $$E(T)=\cup_{i=1}^m \cup_{f\in f(e_i)}(e_i,f)$$
By hypothesis, $|f(e_i)|=d(e_i)=2$ for all $1\le i\le m$.
We then define $\Phi(T)$ to be the graph with vertex set $V_f$ and edge set $$E(\Phi(T))=\cup_{i=1}^m\{(f,f'):f(e_i)=\{f,f'\}\}$$
In other words, $\Phi$ replaces each pair of edges $${(e_i,f), f\in f(e_i)}$$ in the tree $T$ with a single edge $$(f,f'), f(e_i)=\{f,f'\}$$ in the graph $\Phi(T)$.
We see that $\Phi(T)$ has $m+1$ vertices and $m$ edges.  To check it is a tree, we just need to show that it is connected.

We note that for $T\in T_{m,m+1;2}$ and any $f_i,f_j\in V_f$, there exists a unique path of the form $$(f_i,e_{r_{1}}),(e_{r_1},f_{s_1}),...,(f_{s_{k-1}},e_{r_k}),(e_{r_k},f_j)$$ in $T$ since it is a bipartite tree.  Under the map $\Phi$, this path is mapped to $$(f_i,f_{s_1}),\dots,(f_{s_{k-1}},f_j)$$ which is a path in $\Phi(T)$.  So $\Phi(T)$ is a connected graph on $m+1$ vertices with $m$ edges and hence is a tree in $T_{m+1}$.

Note that $S_m\times S_{m+1}$ acts on $T_{m,m+1;2}$ via $S_m$ reordering $V_e$, $S_m+1$ reordering $V_f$ with the corresponding action on edges.  Similarly $S_m\times S_{m+1}$ acts on $T_{m+1}$ via $S_m$ acting trivially, and $S_{m+1}$ reordering $V_f$.  The map $\Phi: T_{m,m+1;2}\to T_{m+1}$ is equivariant with respect to the $S_m\times S_{m+1}$ action.

We observe that a tree $G$ in $T_{m+1}$ has $m$ labeled edges: $\{(f_{s_i},f_{t_i}):1\le i\le m\}$. The preimage of $G$, $\Phi^{-1}(G)$ consists of trees $\{T_{\sigma}:\sigma\in S_m\}$ in $T_{m,m+1;2}$ with edge sets:
$$E(T_{\sigma})=\{(e_{\sigma(i)},f_{s_i}),(e_{\sigma(i)},f_{t_i}):1\le i\le m\}$$
This observation shows that $\Phi:T_{m,m+1;2}\to T_{m+1}$ is a surjective $S_m\times S_{m+1}$-equivariant map. Considering just the $S_m$-action, it shows that the induced map $$T_{m,m+1;2}/S_m\to T_{m+1}/S_m =T_{m+1}$$ is surjective.  But the computation of $\Phi^{-1}(G)$ shows that it consists of one $S_m$-orbit and hence $T_{m,m+1;2}/S_m$ is in bijection with $T_{m+1}$.
Then 
$$T_{m,m+1;2}/(S_{m}\times S_{m+1})=(T_{m,m+1;2}/S_m)/S_{m+1}=T_{m+1}/S_{m+1}=U(m+1)$$
since it is clear that the $S_{m+1}$-orbits of the labeled spanning trees on $m+1$ vertices are in bijection with the set $U(m+1)$ of unlabeled spanning trees on $m+1$ vertices.

\end{proof}

\section{Computing the Cohomology Ring}

An application of Proposition \ref{max cones2} is given by Theorem 12.4.4 of [CLS11]. The theorem, which we state below, gives an explicit way to compute the cohomology ring of a smooth complete toric variety by using primitive collections.

\begin{theorem}[Jurkiewicz-Danilov]\label{cohomology ring}
Let $X_\ss$ be a smooth complete toric variety. For the polynomial ring $\zz [x_1,...,x_r]$ with variables indexed by $\rho_1,...,\rho_r \in \ss (1)$, let $I$ and $J$ be the ideals in $\zz [x_1,...,x_r]$ given by
\begin{align*}
I &= \lb x_{l_1}\cdots x_{l_s} \:|\: l_i \text{ are distinct and } \rho_{l_1} +\cdots+ \rho_{l_s} \text{ is not a cone of } \ss \rb \\
J &=\langle \sum_{l=1}^r \lb m,u_{\rho_l} \rb x_l \:|\: m\in M \rangle
\end{align*}
and define
\begin{align*}
R(\ss)=\zz [x_1,...,x_r]/(I+J).
\end{align*}
Then $x_l\mapsto [D_{\rho_l}]$ induces a ring isomorphism $R(\ss)\simeq H^* (X_{\ss},\zz)$ where $D_{\rho_l}$ is the divisor associated to the ray $\rho_l\in \ss(1)$.
\end{theorem}

The ideal $I$ in the theorem is called the {\it Stanley-Reisner} ideal and, by its definition, one can determine that it is generated by products that correspond to primitive collections. Specifically, the generators of $I$ are the monomials
\begin{align*}
&\prod_{\rho_l\in P} x_l, && P\text{ is a primitive collection for }\ss.
\end{align*}

For $X=X_{mn}$ with $n=1\mod m$, we have that $r$ is equal to the product $mn$ and we can index our variables by $(i,j)\in \llm\times\lln$ where $x_{ij}$ corresponds to the ray given by $e_{ij}$. By Proposition \ref{max cones2}, we have that the primitive collections for $\ss_{mn}$ can be written as
$$P=\{e_{ij} \: |\: i\in Y,j\in Z\}$$
where $Y\subset \{1,\ldots m\}$ and $Z\subset \{1,\ldots ,n\}$ with $1\leq |Y|\leq m$ and $|Z|=(m-|Y|)k+1$. This implies that we can write the generators of $I$ as
\begin{align}
f(Y,Z)&=\prod_{i\in Y}\prod_{j\in Z} x_{ij} \nonumber \\
&=\prod_{(i,j)\in Y\times Z} x_{ij} \label{Igenerators}
\end{align}
for such sets $Y$ and $Z$.

Now that we have a description of the generators for $I$, let us consider the ideal $J$. We can find generators for $J$ by considering any basis for $M$, so let us take the standard basis
$$B=\{e_{ij}\:|\: 1\leq i\leq m-1,1\leq j\leq n-1\}$$ and compute for each $e_{i'j'}\in B$:
\begin{align*}
\sum_{l=1}^r \lb e_{i'j'},u_{\rho_l} \rb x_l&=\sum_{i=1}^{m} \sum_{j=1}^{n} \lb e_{i'j'},e_{ij} \rb x_{ij} \\
&=x_{i'j'}-x_{i'n}-x_{mj'}+x_{mn}.
\end{align*}
This implies that $J$ is the ideal
\begin{align*}
J=\lb x_{ij}-x_{in}-x_{mj}+x_{mn} \:|\: 1\leq i\leq m-1,1\leq j\leq n-1 \rb
\end{align*}
and therefore we have
$$x_{ij}=x_{in}+x_{mj}-x_{mn}$$
in $R(\ss)$ for all $1\leq i\leq m-1$ and $1\leq j\leq n-1$. In particular, if we let
\begin{align*}
S&=\zz[x_{ij}\:|\:1\leq i\leq m,1\leq j\leq n]/J, \\
S'&=\zz[x_{1n},\ldots x_{m-1,n},x_{m,1},\ldots ,x_{m,n-1}, x_{m,n}],
\end{align*}
then there is an isomorphism
\begin{align*}
\phi: S/J&\to S' \\
x_{ij}&\mapsto x_{in}+x_{mj}-x_{mn}
\end{align*}
so we can now compute $R(\ss)=S/(I+J)$ as the quotient $R(\ss)=S'/\phi(I)$.

Using the generators for $I$ given by \eqref{Igenerators}, we have that $\phi(I)$ is generated by polynomials
\begin{align}
\phi(f(Y,Z))&=\prod_{(i,j)\in Y\times Z} \phi(x_{ij}) \nonumber \\
&=\prod_{(i,j)\in Y\times Z} (x_{in}+x_{mj}-x_{mn}) \label{phiIgenerators}
\end{align}
where $Y\subset \llm$ is non-empty and $Z\subset \lln$ is of size $(m-|Y|)k+1$. We claim that we can write $\phi(f(Y,Z))$ as
\begin{align*}
\phi(f(Y,Z))&=\sum_{\substack{A,B\subset Y\times Z \\ A\cap B=\emptyset}} (\prod_{i\in Y} x_{in}^{a_i})(\prod_{j\in Z} x_{mj}^{b_j})(-x_{mn})^{yz-a-b}
\end{align*}
where
\begin{align*}
y&=|Y| \\
z&=|Z| \\
a&=|A| \\
b&=|B| \\
a_i&=|A\cap (\{i\}\times Z)| \\
b_j&=|B\cap (Y\times \{j\})|.
\end{align*}

To see why, consider what a monomial looks like in the expansion of \eqref{phiIgenerators}. Each factor of the monomials come from one of the terms $x_{in},x_{mj},-x_{mn}$ and we get one factor for each $(i,j)\in Y\times Z$. The sets $A$ and $B$ correspond to the subsets of $Y\times Z$ for which the factors $x_{in}$ and $x_{mj}$, respectively, are present. Then the degree of $(-x_{mn})$ is the result of the fact that each monomial must have degree equal to $|Y\times Z|$.

Thus, Theorem \ref{cohomology ring} gives us the following description of the cohomology ring of $X_{mn}$.

\begin{corollary}
The cohomology ring for $X_{mn}$ can be described by the isomorphism
$$H^* (X_{\ss},\zz)\simeq S'/\langle \sum_{\substack{A,B\subset Y\times Z \\ A\cap B=\emptyset}} (\prod_{i\in Y} x_{in}^{a_i})(\prod_{j\in Z} x_{mj}^{b_j})(-x_{mn})^{yz-a-b} \:|\: z=(m-y)k+1 \rangle$$
where the correspondence is given by
\begin{align*}
[D_{e_{ij}}]&\mapsto \begin{cases}
x_{ij} &\text{if } i=m \text{ or } j=n, \\
x_{in}+x_{mj}-x_{mn} &\text{if } i\neq m \text{ and } j\neq n.
\end{cases}
\end{align*}
\end{corollary}

\addcontentsline{toc}{chapter}{Bibliography}

\bibliographystyle{alpha} 

\bibliography{klyachkorefs}{}

\end{document}